\documentclass[12pt]{amsart}
\usepackage{amsmath,amsfonts,euscript,amscd,amsthm,amssymb,upref,graphics,color}

\theoremstyle{plain}
\newtheorem{theorem}{Theorem}
\swapnumbers

\newtheorem{lemma}[subsection]{Lemma}

\theoremstyle{definition}

\newtheorem{nothing*}[subsection]{}

\newcommand{\rien}[1]{}

\newcommand{\C}{\ensuremath{\mathbb{C}}}

\def\e{\epsilon}

\renewcommand{\epsilon}{\varepsilon}
\renewcommand{\phi}{\varphi}
\renewcommand{\emptyset}{\varnothing}
\begin{document}
\renewcommand{\baselinestretch}{1.07}
%%%%%%  TOPMATTER:   %%%%%%%%%%%%%%%%%%%%%%%%%

\title[A Fatou-Bieberbach domain in $\C^2$ which is not Runge]
{A Fatou-Bieberbach domain in $\C^2$ which is not Runge}

\author{Erlend Forn\ae ss Wold}

\address{Matematisk Institutt \\ Universitetet i Oslo
\\Postboks 1053 Blindern
\\ NO-0316 Oslo, Norway}
\email{erlendfw@math.uio.no}

\date{January 14, 2007}
\thanks{Supported by Schweizerische Nationalfonds grant 200021-116165/1}

 \begin{abstract}
 Since a paper by J.P.Rosay and W.Rudin from 1988 there has been an open question
 whether all Fatou-Bieberbach domains are Runge.  We give an example of a
 Fatou-Bieberbach domain $\Omega$ in $\C^2$ which is not Runge.  The domain $\Omega$
 provides (yet) a negative answer to a problem of Bremermann.

\end{abstract}
\maketitle \vfuzz=2pt

\vfuzz=2pt
%%%%%%%%%%%%%%%%%%%%%%%%%%%%%%%%%%%%%%%%%%%%%%%%%%%%%%%%%%%%%%%%%%%
%%%%%%%%%%%%%%%%%%%%%%%%%%%%%%%%%%%%%%%%%%%%%%%%%%%%%%%%%%%%%%%%%%%
%%%%%%%%%%%%%%%%%%%%%%%%%%%%%%%%%%%%%%%%%%%%%%%%%%%%%%%%%%%%%%%%%%%

\vskip 1cm

\section{Introduction}

\bigskip

We give a negative answer to the problem, initially posed by J.P.
Rosay and W. Rudin in \cite{RR} and later in \cite{FR}, as to
whether all Fatou-Bieberbach domains are Runge:

\bigskip

\begin{theorem}\label{main}
There is a Fatou-Bieberbach domain $\Omega$ in $\C^*\times\C$ which
is Runge in $\C^*\times\C$ but not in $\C^2$.
\end{theorem}

\bigskip

A Fatou-Bieberbach domain is a proper subdomain of $\C^n$ which is
biholomorphic to $\C^n$, and a domain $\Omega\subset\C^n$ is said to
be Runge (in $\C^n$) if any holomorphic function $f\in\mathcal{O}(\Omega)$ can
be approximated uniformly on compacts in $\Omega$  by polynomials. \

It should be noted that although the domain $\Omega$ is not Runge it
still has the property that the intersection of $\Omega$ with any
complex line $L$ is simply connected: Let $V$ be a connected
component of $\Omega\cap L$, let $\Gamma\subset V$ be a simple
closed curve , and let $D$ denote the disk in $L$ bounded by
$\Gamma$.  Since $\Gamma$ is null-homotopic in $\Omega$ we have that
$D$ is contained in $\C^*\times\C$ and so the claim follows from the
fact that $\Omega$ is Runge in $\C^*\times\C$.  Intersecting
$\Omega$ with a suitable bounded subset of $\C^2$ this gives a
negative answer to the problem of Bremermann: "Suppose that $D$ is a
Stein domain in $\mathbb C^n$ such that for every complex line $l$
in $\mathbb{C}^n$, $l\setminus D$ is connected. Is it true that $D$
is Runge in $\mathbb C^n$?". Negative answers to this problem have
also recently been given in \cite{Ab} and \cite{Jo}. On can infact
show, using an argument as above together with the argument
principle, that if $\mathcal{R}$ is a smoothly bounded planar domain
and if $\phi(\mathcal{R})$  is a holomorphic embedding of
$\mathcal{R}$ into $\C^2$  with
$\phi(\partial\mathcal{R})\subset\Omega$, then
$\phi(\mathcal{R})\subset\Omega$.

\medskip

The idea of the proof is the following: Observe first that if
$\Omega$ is a Fatou-Bieberbach domain in $\C^2$ which is Runge,
then $\Omega$ has the property that if $Y\subset\Omega$ is compact
then its polynomially convex hull
$$
\widehat Y:=\{(z,w)\in\C^2;|P(z,w)|\leq\|P\|_Y \ \forall
P\in\mathcal{P}(\C^2)\}
$$
is contained in $\Omega$.  To prove the theorem we will construct a
domain $\Omega$  such that $\widehat Y\setminus\Omega\neq\emptyset$
for a certain compact set $Y$. For a compact subset
$Y\subset\C^*\times\C$  let $\widehat Y_*$  denote the set
$$
\widehat Y_*:=\{(z,w)\in\C^2;|P(z,w)|\leq\|P\|_Y \ \forall
P\in\mathcal{O}(\C^*\times\C)\}.
$$
We say that the set $Y$ is holomorphically convex if $\widehat
Y_*=Y$.  We will first construct (a construction by Stolzenberg) a
holomorphically convex compact set $Y\subset\C^*\times\C$ having the
property that $\widehat Y\cap(\{0\}\times\C)\neq\emptyset$. $Y$ is
the disjoint union of two disks is $\C^*\times\C$.  We will then use
the fact that $\C^*\times\C$ has the \emph{density property} to
construct a Fatou-Bieberbach domain $\Omega\subset\C^*\times\C$ such
that $Y\subset\Omega$. The domain $\Omega$  cannot be Runge.

\medskip

A few words about the density property and approximation by
automorphisms. As defined in \cite{Va}, \emph{a complex manifold $M$
is said to have the density property if every holomorphic vector
field on $M$ can be approximated locally uniformly by Lie
combinations of complete vector fields on $M$}.  It was proved in
\cite{Va} that $\C^*\times\C$  has the density property. In
Anders\'{e}n-Lempert theory the density property corresponds to the
fact that in $\C^n$ every entire vector field can be approximated by
sums of complete vector fields. This has been studied also in
\cite{Va2}. \

Using the density property of $\C^*\times\C$ one gets as in
\cite{FR} (by copying their arguments): Let $\Omega$ be an open set
in $\C^*\times\C$. For every $t\in [0,1]$, let $\phi_t$ be a
biholomorphic map from $\Omega$ into $\C^*\times\C$, of class
$\mathcal{C}^2$ in $(t,z)\in [0,1]\times\Omega$.  Assume that
$\phi_0=Id$, and assume that each domain $\Omega_t=\phi_t(\Omega)$
is Runge in $\C^*\times\C$.  Then for every $t\in [0,1]$ the map
$\phi_t$ can be approximated on $\Omega$ by holomorphic
automorphisms of $\C^*\times\C$. In the proof of Theorem 1 we will
construct such an isotopy. \

We will let $\pi$ denote the projection onto the first coordinate in
$\C^*\times\C$  and in $\C^2$, and we will let $B_\e(p)$ denote the
open ball of radius $\e$ centered at a point $p$. \

\bigskip

\section{construction of the set Y}

\bigskip

We start by defining a certain rationally convex subset $Y$ of
$\C^2$. The set will be a union of two disjoint polynomially convex
disks in $\C^*\times\C$, but the polynomial hull of the union will
contain the origin. This construction is taken from \cite{St}, page
392-396, and is due to Stolzenberg \cite{Sb}. \

Let $\Omega_1$ and $\Omega_2$ be simply connected domains in $\C$,
as in Fig.1. below, with smooth boundary, such that if
$I_+=[1,\sqrt{3}],I_-=[-\sqrt{3},-1]$, then
$I_+\subset\partial\Omega_1$, $I_-\subset\partial\Omega_2$. Require
that $\partial\Omega_1$  and $\partial\Omega_2$  meet only twice,
that $I_-\subset\Omega_1, I_+\subset\Omega_2$, and, finally, that
$\partial\Omega_1\cup\partial\Omega_2$ be the union of the boundary
of the unbounded component of
$\C\setminus(\partial\Omega_1\cup\partial\Omega_2)$, together with
the boundary of the component of this set that contains the origin.
Let the intersections of the boundaries be the points $i$  and $-i$.
\

%TeXCAD Picture [bilde3.pic]. Options:
%\grade{\on}
%\emlines{\off}
%\epic{\off}
%\beziermacro{\on}
%\reduce{\on}
%\snapping{\off}
%\pvinsert{% Your \input, \def, etc. here}
%\quality{8.00}
%\graddiff{0.01}
%\snapasp{1}
%\zoom{4.0000}
\unitlength 1mm % = 2.85pt
\linethickness{0.4pt}
\ifx\plotpoint\undefined\newsavebox{\plotpoint}\fi % GNUPLOT compatibility
\begin{picture}(112.25,120.5)(-10,0)
\put(37,83){\line(1,0){13.5}} \put(81,83){\line(-15,0){13.5}}
\qbezier(37,83)(11.5,83.38)(60,94.25)
\qbezier(81,83)(106.5,83.38)(58,94.25)
\qbezier(60,94)(99,102.38)(93,82.25)
\qbezier(58,94)(19,102.38)(25,82.25)
\qbezier(50.25,83)(55.88,84.63)(60,66.75)
\qbezier(67.75,83)(62.12,84.63)(58,66.75)
\qbezier(60,66.75)(72,14.5)(93,82.25)
\qbezier(58,66.75)(46,14.5)(25,82.25) \thicklines
\put(12,83.25){\line(1,0){94.25}} \put(59,48.25){\line(0,1){69.75}}
\put(112.25,83){\makebox(0,0)[cc]{$Rez$}}
\put(59,120.5){\makebox(0,0)[cc]{$Imz$}}
\put(60.5,97.25){\makebox(0,0)[cc]{$i$}}
\put(61,70.75){\makebox(0,0)[cc]{$-i$}}
\put(51,82.25){\line(0,1){2}} \put(31.25,82.25){\line(0,1){2}}
\put(67,82){\line(0,1){2}} \put(86.5,82.25){\line(0,1){2}}
\put(31.5,79){\makebox(0,0)[cc]{$-\sqrt{3}$}}
\put(67.75,79){\makebox(0,0)[cc]{$1$}}
\put(50.75,79.25){\makebox(0,0)[cc]{$-1$}}
\put(87,79.25){\makebox(0,0)[cc]{$\sqrt{3}$}}
\put(77.25,85.25){\makebox(0,0)[cc]{$I_+$}}
\put(42.25,85.25){\makebox(0,0)[cc]{$I_-$}}
\put(72.75,61.75){\makebox(0,0)[cc]{$\Omega_2$}}
\put(45.25,61.75){\makebox(0,0)[cc]{$\Omega_1$}}
\put(97,77.75){\makebox(0,0)[cc]{$\partial\Omega_2$}}
\put(21.5,79.5){\makebox(0,0)[cc]{$\partial\Omega_1$}}
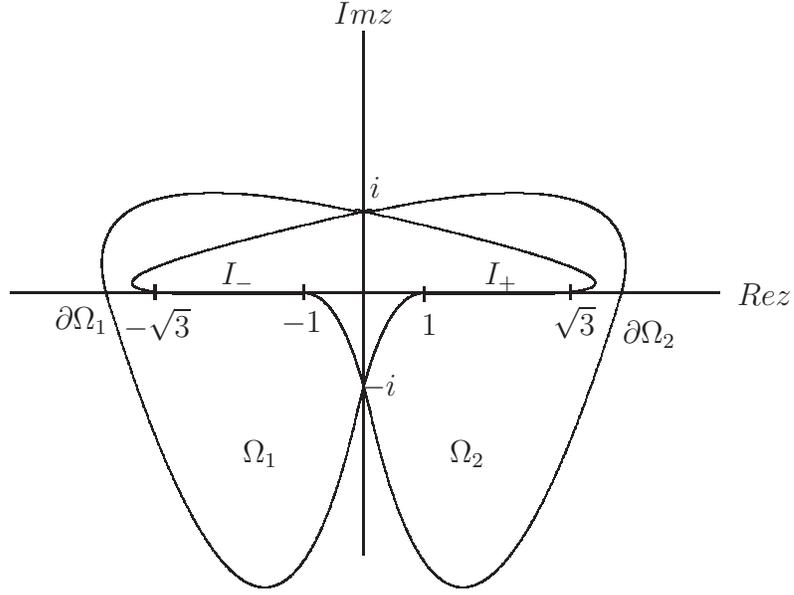
\put(58.00,40){\makebox(0,0)[cc]{{\tiny\textsc{Fig.1.}Two smoothly
bounded simply connected domains $\Omega_j$ with
$\partial\Omega_1\cap\partial\Omega_2=\pm i$.}}}
\end{picture}

We define

\

\centerline{$V_1=\{(z,w)\in\C^2; z^2-w \ \mathrm{is \ real \ and \
lies \ in} \ [0,1]\}$,} \

\centerline{$V_2=\{(z,w)\in\C^2; w \ \mathrm{is \ real \ and \ lies
\ in} \ [1,2]\}$,} \

\centerline{$X_1=\{(z,w)\in V_1;z\in\partial\Omega_2\}$,} \

\centerline{$X_2=\{(z,w)\in V_2;z\in\partial\Omega_1\}$,} \

\

Note that $X_1$  and $X_2$  are totally real annuli, that they are
disjoint, and that the origin is contained in the polynomial hull of
$X_1$.  Next we want to remove pieces from $X_1$  and $X_2$ to
create two disks. \

Define

\

\centerline{$\tilde V_1=V_1\cap\pi^{-1}(I_+),$}\

\centerline{$\tilde V_2=V_2\cap\pi^{-1}(I_-),$} \

\centerline{$Y_1=\overline{X_1\setminus\tilde V_2}$,} \

\centerline{$Y_2=\overline{X_2\setminus\tilde V_1}$.}

\

The set $Y$  will be defined as $Y=Y_1\cup Y_2$. Note that
$$
(*) \ \tilde V_1\subset\widehat X_1, \tilde V_2\subset\widehat X_2.
$$
Let us describe what $Y_1$  and $Y_2$  looks like over $I_-$ and
$I_+$ respectively.  By the equations we see that these sets are
contained in $\mathbb{R}^2$.  Let $(x,y)$ denote the real parts of
$(z,w)$. \

Over $I_-$ we have that $Y_1$ is the union of the two sets defined
by

\

$(a)$ \centerline{$2\leq y\leq x^2 \ \mathrm{if} \ -\sqrt{3}\leq
x\leq-\sqrt{2},$} \

$(b)$ \centerline{$x^2-1\leq y\leq 1 \ \mathrm{if} \ -\sqrt{2}\leq
x\leq-1.$}

\

Over $I_+$ we have that $Y_2$ is the union of the sets defined by

\

$(c)$ \centerline{$x^2\leq y\leq 2 \ \mathrm{if} \ 1\leq
x\leq\sqrt{2},$} \

$(d)$ \centerline{$1\leq y\leq x^2-1 \ \mathrm{if} \ \sqrt{2}\leq
x\leq\sqrt{3}.$}

\

From these equations we see that $Y_1$  and $Y_2$ are disks.

We have that
$$
(**) \ \widehat Y \ \mathrm{contains \ the \ origin}
$$
because of the following: We already noted that the origin is
contained in $\widehat X_1$, so the claim follows from $(*)$ and the
following simpler version of Lemma 29.31, \cite{St}, page 392: Let $X_1$
and $X_2$ be disjoint compact sets in $\C^N$, and let $S_1$ and
$S_2$ be relatively open subsets of $X_1$ and $X_2$ respectively
such that $S_1\subset\widehat X_2, S_2\subset\widehat X_1$. Then
$\widehat{X_1\cup X_2}=\widehat{(X_1\setminus S_1)\cup (X_2\setminus
S_2)}$. The reason for this, which was pointed out by the referee,
is simply that neither $S_1$ nor $S_2$ can contain peak points for
the algebra generated by the polynomials on $X_1\cup X_2$. \

\section{Proof of Theorem \ref{main}}

It is proved in \cite{St} that the set $Y$ is rationally convex, and
that the sets $Y_j$ are polynomially convex separately. For our
construction we need to know that $Y$ is holomorphically convex, so
we prove the following:

\begin{lemma}\label{hc}
We have that $Y$ is holomorphically convex in $\C^*\times\C$.
\end{lemma}
\begin{proof}
For $j=1,2,$ let $Y^+_j$  and $Y^-_j$ denote the sets
$Y_j\cap\{\mathrm{Re}(z)\geq 0\}$ and $Y_j\cap\{\mathrm{Re}(z)\leq
0\}$ respectively.  Let $Y^+=Y^+_1\cup Y^+_2$ and $Y^-=Y^-_1\cup
Y^-_2$. \

Observe first that $Y^+$ and $Y^-$ are polynomially convex
separately: Assume to get a contradiction that $\widehat{Y^-}$
contains nontrivial points. In that case there exists a graph $G(f)$
of a bounded holomorphic function defined on the topological disk
$U$ bounded by $\pi(Y^-)$, such that $G(f)\subset \widehat{Y^-}$,
and such that $(z,f(z))\in Y^-$  for a.a. (in terms of radial limits
if we regard $U$ as a proper disk) $z\in\pi(Y^-)$ (Theorem 20.2. in
\cite{AW}, page 172, holds by the discussion on page 171 even though
the fibers over $\pm i$ are not convex). Then for continuity reasons
$\overline{G(f)}$ would have to contain nontrivial points of
$\widehat{Y^-}$ in the fibers $\{\pm i\}\times\mathbb C$ - but as
this clearly cannot be the case, we have our contradiction. The case
of $Y^+$ is similar. \

Next assume to get a contradiction that there is a point
$(z_0,w_0)\in\widehat Y_*\setminus Y$ with $\mathrm{Re}(z_0)<0$. The
function $f(z)$ defined to be $(z+i)(z-i)$ on $\pi(Y^-)\cup\{z_0\}$
and zero on $\pi(Y^+)$ can be uniformly approximated on
$\pi(Y)\cup\{z_0\}$ by polynomials in $z$ and $\frac{1}{z}$, and so
any representing Jensen measure (see \cite{St} Chapter 2) for the
functional $g\mapsto g(z_0,w_0)$ would have to be supported on
$Y^-$.  But then the point $(z_0,w_0)$ would have to be in the hull
of $Y^-$ which is a contradiction. The corresponding conclusion
holds for $\mathrm{Re}(z_0)>0$. \

Finally, Rossi's local maximum principle excludes the possibility of
there being nontrivial points in the hull contained in $\{\pm
i\}\times\mathbb C$.
\end{proof}

\begin{lemma}\label{aut}
Let $p=(z_0,w_0)\in\C^*\times\C$  and let $\e>0$.  Then there
exists an automorphism $\psi$ of $\C^*\times\C$  such that
$\psi(Y)\subset B_\e(p)$.
\end{lemma}
\begin{proof}

We need to argue that there exists an isotopy as described in the
introduction, and we content ourselves by demonstrating that there
exist isotopies mapping $Y_1$ and $Y_2$ into separate arbitrarily
small balls - the rest is trivial.  Let $q_j\in Y_j$ be a point for
$j=1,2,$ and let $\delta>0$.  Since $Y_j$ is a smooth disk there
clearly exists a smooth map $f^j:[0,1]\times Y_j\rightarrow Y_j$
such that for each fixed $t$ the map $f^j_t:Y_j\rightarrow Y_j$ is a
smooth diffeomorphism, such that $f^j_0$ is the identity, and such
that $f^j_1(Y_j)\subset B_\delta(q_j)$.  Since $Y_j$ is totally real
there exists, by \cite{fo} Corollary 3.2, for each $\epsilon>0$ a
real analytic map $\Phi^j:[0,1]\times\mathbb C^2\rightarrow \mathbb
C^2 $ such that $\Phi^j_t\in\mathrm{Aut_{hol}}(\C^2)$ for each $t$,
$\Phi^j_0$ is the identity, and $\|f^j-\Phi^j\|_{[0,1]\times
Y_j}<\epsilon$. For small enough $\epsilon$ we restrict $\Phi^j$ to
a sufficiently small Runge neighborhood of $Y_j$.

\end{proof}

\emph{Proof of Theorem 1}: Let $G$ be an automorphism of
$\C^*\times\C$ with an attracting fixed point $p\in\C^*\times\C$. It
is well known that the basin of attraction of the point $p$ is a
Fatou-Bieberbach domain.  This domain is clearly contained in
$\C^*\times\C$.  Denote this domain by $\Omega(G)$.  Let $\e$ be a
positive real number such that $B_\e(p)\subset\Omega(G)$.  By Lemma
\ref{aut} there is an automorphism $\psi$ of $\C^*\times\C$ such
that $\psi(Y)\subset B_\e(p)$.  Then $Y\subset\psi^{-1}(\Omega(G))$.
The set $\psi^{-1}(\Omega(G))$ is biholomorphic to $\C^2$, and from
$(**)$ in Section 2 we have that $\widehat Y$  contains the origin.
On the other hand it is clear that $\Omega(G)$ is Runge in
$\C^*\times\C$, and so $\psi^{-1}(\Omega(G))$ is Runge in
$\C^*\times\C$.

$\hfill\square$

\medskip

\emph{Acknowlegements.}  The author would like to thank the referee for comments; in particular
for suggesting a simplified proof of Lemma 3.1.  

\bibliographystyle{amsplain}

 \end{document}